# RECURSIVE PARTITION STRUCTURES

By Alexander V. Gnedin and Yuri Yakubovich[1]

*Utrecht University*

A class of random discrete distributions $P$ is introduced by means of a recursive splitting of unity. Assuming supercritical branching, we show that for partitions induced by sampling from such $P$ a power growth of the number of blocks is typical. Some known and some new partition structures appear when $P$ is induced by a Dirichlet splitting.

**1. Introduction.** By a random discrete distribution (or a *paintbox*) we shall understand an infinite collection $P = (P_j)$ of nonnegative random variables whose sum is unity. Interpreting the terms of $P$ as frequencies of distinct colors, Kingman's paintbox construction [16] defines a random exchangeable partition $\mathcal{P}$ of an infinite set of balls labeled $1, 2, \ldots$ in such a way that, conditionally given $(P_j)$, the generic ball $n$ is painted color $j$ with probability $P_j$, independently of all other balls. The blocks of $\mathcal{P}$ are composed of balls painted the same color. Two paintboxes which only differ by the arrangement of terms in a sequence yield the same $\mathcal{P}$; hence to maintain symmetry we may identify the paintboxes with the point process $\sum_j \delta_{P_j}$. See [3, 22] for extensive background on exchangeable partitions.

Let $K_{nr}$ be the number of colors represented exactly $r$ times on $n$ first balls, and let $K_n$ be the total number of different colors represented on $n$ first balls, so that $\sum_r K_{nr} = K_n$, $\sum_r r K_{nr} = n$. The sequence of joint distributions of $(K_{n1}, \ldots, K_{nn})$ for $n = 1, 2, \ldots$ is a *partition structure*, that is, a consistent family of distributions on partitions [16]. We are interested in the asymptotic features of $K_n$ and $K_{nr}$'s, as $n \to \infty$, for one particular class of models for $P$.

The functionals $K_n$ and $K_{nr}$ have been studied in some depth for several families of random discrete distributions. The best known instance is the

Received October 2005; revised February 2006.
[1]Supported by the Netherlands Organization for Scientific Research (NWO).
*AMS 2000 subject classifications.* 60G09, 60C05.
*Key words and phrases.* Partition structure, recursive construction, Dirichlet distribution.

Poisson–Dirichlet/GEM paintbox, which induces $\mathcal{P}$ called the Ewens partition. For the Ewens partition $K_n$ is approximately Gaussian with the mean and the variance both growing logarithmically, and the sequence $(K_{n1}, K_{n2}, \ldots)$ converges, as $n \to \infty$, to a sequence of independent Poisson variates [1]. The GEM realization of the Poisson–Dirichlet paintbox amounts to the *stick-breaking* representation $P_j = W_1 \cdots W_{j-1}(1 - W_j)$ $(j = 1, 2, \ldots)$ with independent $W_j$'s distributed according to beta$(\theta, 1)$. A larger class of models for $P$ of this type, with arbitrary independent identically distributed factors $W_j \in [0, 1]$, was studied in [6], where it was shown that, under very mild assumptions on the distribution of $W_j$'s, the behavior of $K_n$ is analogous to that in the Ewens case.

Each random discrete distribution $P$ resulting from the stick-breaking can be viewed as a collection of jump sizes of the process $(\exp(-S_t), t \geq 0)$, where $(S_t)$ is a compound Poisson process. A considerable extension of this scheme (see [7, 8]) appears when we assume $(S_t)$ to be a subordinator with some infinite Lévy measure $\nu$. It is known that the orders of growth of $K_n$ and $K_{nr}$'s are determined then by the behavior of the tail $\nu[x, \infty[$ for $x \downarrow 0$. Specifically, if the tail behaves like $x^\alpha \ell(1/x)$ with $0 < \alpha < 1$ and $\ell$ a function of slow variation at infinity, then the order of growth of $K_n$ and all $K_{nr}$'s is $n^\alpha \ell(n)$, and with this scaling $K_n$ and each $K_{nr}$ converge, almost surely, to constant multiples of the same random variable [10]. A distinguished example of the latter situation is the Poisson–Dirichlet paintbox [23] with two parameters $0 < \alpha < 1$ and $\theta > -\alpha$, which induces the Ewens–Pitman partition structure whose distribution is given by the formula

$$\mathbb{P}[K_{n1} = k_1, \ldots, K_{nn} = k_n]$$
$$= n! \frac{\theta(\theta + \alpha) \ldots (\theta + (\ell - 1)\alpha)}{(\theta)_{n\uparrow}} \prod_{i=1}^n \frac{((1 - \alpha)_{i-1})^{k_i}}{i!^{k_i} k_i!},$$

where $(a)_{i\uparrow} = a(a + 1) \ldots (a + i - 1)$ stands for rising factorials, and $\ell = k_1 + \cdots + k_n$. A construction of this partition structure via $\exp(-S_t)$ is given in [7], and in Section 6.1 we briefly recall the original construction [20] of the two-parameter paintbox. Very different asymptotic behavior appears when the tail $\nu[x, \infty[$ is slowly varying at zero like, for example, for gamma subordinators: in this case the moments of $K_n$ and $K_{nr}$'s are slowly varying functions of $n$, all $K_{nr}$'s grow on the same scale but slower than $K_n$ and, subject to a suitable normalization, $K_n$ is asymptotically Gaussian [2, 11].

The stick-breaking model for $P$ is the simplest instance of a recursive construction in the sense of the present paper. By stick-breaking the unity splits in two pieces $1 - W_1$ and $W_1$, the first piece becomes a term of $P$, and the second keeps on dividing by the same rule, hence producing again a term for $P$ and a piece which divides further, and so on. The paintbox associated with $\exp(-S_t)$, for $(S_t)$ a subordinator with infinite Lévy measure, can be also realized via a recursive construction which produces, at each step, infinitely many terms for $P$ and exactly one piece to undergo further splitting.

In this paper, the principal step away from the models of the stick-breaking type is that we deal with the recursive models for $P$ in which a random splitting of unity involves some branching. We assume an inductive procedure in which each step yields a collection of terms included in $P$, and a multitude of divisible pieces to iterate a random splitting rule. A typical example of our class of models is the paintbox arising by the following construction of a random Cantor set $\mathcal{R}$ (see [18] for the general theory of recursive constructions of this kind). Start with dividing the unit interval in three intervals of sizes $X_1, Y_2, X_3$, from the left to the right, as obtained by cutting $[0, 1]$ at the locations of two uniform order statistics. Remove the middle interval and iterate the operation of cutting and removing the middle independently on two other intervals (considered as scaled copies of $[0, 1]$), then iterate on four intervals, and so on. A random set $\mathcal{R}$ of Lebesgue measure zero is defined as the complement to the union of all removed intervals, and the collection of lengths of the removed intervals arranged in some sequence defines a paintbox $P$. It will follow from the main result of this paper (Theorem 5) that $K_n$ and each $K_{nr}$ grow for this $P$ like $n^{\alpha_*}$ with exponent $\alpha_* = (\sqrt{17} - 3)/2$, which is equal to the Malthusian parameter of a related branching process and is also equal to the Hausdorff dimension of $\mathcal{R}$ [18, 19].

In wider terms, our construction is described as follows. At step one the unity is randomly divided in some collection of *solids* and some collection of *crumbs*. The solids immediately suspend further transformation while the crumbs keep on falling apart. At step two the crumbs are split further by the same random rule, the newly created solids become indivisible and the crumbs are subject to further division, and so on. Eventually, the crumbs decompose completely in solids, and the sizes of solids (arranged in some sequence) comprise the paintbox $P$.

We will show that the power growth of $K_n$ and $K_{nr}$'s is quite common for $\mathcal{P}$ derived from such a recursive paintbox with supercritical branching. Moreover, by the power scaling $K_n$ and $K_{nr}$'s all converge to constant multiples of the same random variable $M$ which can be characterized in terms of a distributional fixed-point equation. Some explicit moment computations for $M$ are possible for instances of the splitting procedure based on the Dirichlet distribution; for some choices of the parameters these yield the paintboxes of the Poisson–Dirichlet $(\alpha, \alpha/d)$-type $(d = 1, 2, \ldots)$ and for some other choices yield new paintboxes hence novel partition structures.

**2. Malthusian hypothesis and a martingale.** Let $(\mathbf{X}, \mathbf{Y}) = ((X_i), (Y_j))$ be two sequences of random variables with values in $[0, 1]$. To introduce a genealogical structure of the division process it will be convenient to assume that the sequences are labeled by two disjoint subsets of $\mathbb{N}$. We further require that

$$\text{(1)} \quad \sum_i X_i + \sum_j Y_j = 1, \qquad \mathbb{E}\left[\sum_i X_i\right] < 1, \qquad \mathbb{E}[\#\{i : X_i > 0\}] > 1.$$

The division process starts with a sole unit crumb $\xi_\varnothing$ (generation 0) which produces the first generation of crumbs $(\xi_i)$ and solids $(\eta_j)$ whose joint law and labels are the same as for $(\mathbf{X}, \mathbf{Y})$. Inductively, the offspring of generation $k-1$ are crumbs $(\xi_{i_1,\ldots,i_{k-1},i})$ and solids $(\eta_{i_1,\ldots,i_{k-1},j})$. The solids stop division, while each crumb $\xi_{i_1,\ldots,i_k}$ splits further into crumbs $(\xi_{i_1,\ldots,i_k,i})$ and solids $(\eta_{i_1,\ldots,i_k,j})$ whose labeling and sizes relative to the parent crumb follow the law of $(\mathbf{X}, \mathbf{Y})$, independently of the history and the sizes of other members of the current generation. The first two assumptions in (1) guarantee that the total size of solids over all generations is unity, hence these sizes (arranged in a sequence) define a paintbox $P$. The third assumption in (1) says that the branching of crumbs is *supercritical*.

Introduce the intensity measures $\sigma$ and $\nu$ by requiring the equalities

$$\mathbb{E}\left[\sum_i f(X_i)\right] = \int_0^1 f(x)\sigma(dx),$$

$$\mathbb{E}\left[\sum_j f(Y_j)\right] = \int_0^1 f(x)\nu(dx)$$

to hold for all nonnegative measurable functions $f$. Substituting power functions $f(x) = x^\alpha$ in these formulas yields the Mellin transforms of the measures

$$\psi(\alpha) := \int_0^1 x^\alpha \sigma(dx),$$

$$\varphi(\alpha) := \int_0^1 x^\alpha \nu(dx).$$

Recall that, as a function of complex parameter, the Mellin transform of a measure on $[0,1]$ is analytical in the half-plane to the right of the convergence abscissa, has a ridge on $\operatorname{Im}\alpha = 0$ and decreases on the real half-line.

The *Malthusian hypothesis* accepted in this paper amounts to the assumptions that:

- there exists a solution $\alpha_*$ to the equation

$$\text{(2)} \qquad \psi(\alpha) = 1$$

(which satisfies then $\alpha_* \in \,]0,1[$ since $\psi(1) < 1 < \psi(0)$ by (1)),

- there exists $\varepsilon > 0$ such that $\alpha_*$ is a unique solution to (2) in the half-plane $\{\alpha : \operatorname{Re}\alpha > \alpha_* - \varepsilon\}$ and $\varphi(\alpha_* - \varepsilon) < \infty$.

Following the established tradition in the theory of branching processes we call $\alpha_*$ the *Malthusian exponent*. Obvious sufficient conditions for the Malthusian hypothesis are $\psi(0) < \infty$ and $\varphi(0) < \infty$. Note also that the second part of the Malthusian hypothesis implies that $\sigma$ is not supported by a geometric progression, since otherwise (2) would have infinitely many periodically spaced roots on the line $\operatorname{Re}\alpha = \alpha_*$.

Summing the $\alpha_*$th powers of crumbs in a given generation yields a remarkable process called the *intrinsic martingale* [13]

$$M_k := \sum_{i_1,\ldots,i_k} \xi_{i_1,\ldots,i_k}^{\alpha_*}, \qquad k = 1, 2, \ldots,$$

which, under the Malthusian hypothesis, converges to a terminal value

$$M := \lim_{k \to \infty} M_k$$

with $\mathbb{E}[M] = 1$; see [17]. The limit variable satisfies the distributional fixed-point equation

$$(3) \qquad M \stackrel{d}{=} \sum_i X_i^{\alpha_*} M^{(i)},$$

where $M^{(i)}$ are independent copies of $M$, independent of $\mathbf{X}$. It is known that (3) along with $\mathbb{E}[M] = 1$ uniquely characterizes $M$ [12], Proposition 3(a).

**3. The mean values of counts.** Consider the powered sums of sizes of all solids that make up the paintbox

$$G_\alpha := \sum_j P_j^\alpha = \sum_{k=1}^\infty \sum_{i_1,\ldots,i_{k-1},j} \eta_{i_1,\ldots,i_{k-1},j}^\alpha$$

and let

$$p(\alpha) := \mathbb{E}[G_\alpha].$$

For integer arguments the value $p(n)$ is the probability that $n$ balls are painted the same color. The first-split decomposition of the division process yields the distributional equation

$$(4) \qquad G_\alpha \stackrel{d}{=} \sum_i X_i^\alpha G_\alpha^{(i)} + \sum_j Y_j^\alpha,$$

where $G_\alpha^{(i)}$ are independent copies of $G_\alpha$ which are also independent of $(\mathbf{X}, \mathbf{Y})$. Taking the expectations this implies

$$(5) \qquad p(\alpha) = \frac{\varphi(\alpha)}{1 - \psi(\alpha)} \qquad \text{for } \operatorname{Re}\alpha > \alpha_* - \varepsilon.$$

By the Malthusian hypothesis the expectations involved are finite, and the function $p$ is meromorphic in the half-plane $\operatorname{Re}\alpha > \alpha_* - \varepsilon$, with a sole simple pole at $\alpha_*$. These analytic properties of $p$ provide a background for establishing the growth properties for the mean values of counts $K_n$ and $K_{nr}$.

Conditionally given $(P_j)$ the probability that at least one of $n$ balls is painted color $j$ is $1 - (1 - P_j)^n$, hence recalling the definition of $p$

$$\mathbb{E}[K_n] = \mathbb{E}\left[\sum_j (1 - (1 - P_j)^n)\right] = \sum_{m=1}^n \binom{n}{m}(-1)^{m+1} p(m). \tag{6}$$

In a similar way, computing the chance that color $j$ is represented exactly $r$ times in $n$ balls:

$$\mathbb{E}[K_{nr}] = \binom{n}{r}\mathbb{E}\left[\sum_j (P_j^r (1 - P_j)^{n-r})\right]$$

$$= \binom{n}{r}\sum_{m=0}^{n-r} \binom{n-r}{m}(-1)^m p(m+r). \tag{7}$$

THEOREM 1. *Under assumption* (1) *and the Malthusian hypothesis, the following asymptotics hold:*

$$\mathbb{E}[K_n] = n^{\alpha_*}\frac{\Gamma(-\alpha_*)\varphi(\alpha_*)}{\psi'(\alpha_*)} + O(n^{\alpha_* - \varepsilon}) \quad \text{as } n \to \infty, \tag{8}$$

$$\mathbb{E}[K_{nr}] = n^{\alpha_*}\frac{\Gamma(r - \alpha_*)\varphi(\alpha_*)}{-r!\psi'(\alpha_*)} + O(n^{\alpha_* - \varepsilon}) \quad \text{as } n \to \infty. \tag{9}$$

PROOF. Because the function $p$ is bounded in the half-plane $\operatorname{Re}\alpha > \alpha_* - \varepsilon$, outside any neighborhood of $\alpha_*$, and because $p$ has a simple pole at $\alpha_*$, we can apply the Rice method (see [5], Theorem 2(ii)) to the alternating sum (6) to obtain

$$\sum_{m=1}^n \binom{n}{m}(-1)^{m+1} p(m)$$

$$= \operatorname*{Res}_{\alpha = \alpha_*} p(\alpha)\frac{\Gamma(1 - \alpha)\Gamma(n + 1)}{\alpha \Gamma(n + 1 - \alpha)} + O(n^{\alpha_* - \varepsilon}). \tag{10}$$

The residue at $\alpha_*$ is equal to $-\varphi(\alpha_*)/\psi'(\alpha_*)$, which taken together with $\Gamma(n + a)/\Gamma(n) \sim n^a$ readily yields (8). The result for $K_{nr}$ can be obtained in the same way. Alternatively, observe that the sum in (7) is asymptotic to a constant multiple of the $r$th derivative of (6) in the variable $n$, hence (9) follows from (8) by a Tauberian argument. $\square$

The first-split decomposition shows that the number of colors $K_n$ satisfies a divide-and-conquer recurrence of the form

$$K_n \stackrel{d}{=} \sum_i K_{A_{ni}}^{(i)} + B_n,$$

where $(K_n^{(i)}, n = 1, 2, \ldots)$ are independent copies of $(K_n)$, and the joint law of $(A_{ni}, B_n)_{i \geq 1}$ follows by considering the partition of $n$ induced by the joint paintbox $(X_1, X_2, \ldots; (\sum_j Y_j))$. In other words, given $(\mathbf{X}, \mathbf{Y})$ each of $n$ balls is painted color $i$ with probability $X_i$ and left uncolored with probability $\sum_j Y_j$ then $A_{ni}$ is the number of balls painted color $i$ and $B_n$ is the number of uncolored balls. Because (3) is a limit analogue of this equation, the contraction method [25] can be exploited to show weak convergence of scaled $K_n$. To argue the strong convergence we will apply an indirect approach (also used in [10]) which relates the growth properties of $K_n, K_{nr}$ with the sizes of solids. Let

$$N_x := \#\{j : P_j \geq x\}$$

be the number of solids with size at least $x$.

LEMMA 2. *If the paintbox satisfies $N_x \sim L x^{-\alpha}$ a.s. as $x \downarrow 0$, with $0 < \alpha < 1$ and $L$ a nonnegative random variable, then for $n \to \infty$*

$$K_n / n^\alpha \to \Gamma(1-\alpha) L \qquad a.s.$$

*and*

$$K_{nr}/n^\alpha \to (\alpha \Gamma(r-\alpha)/r!) L \qquad a.s.$$

PROOF. Conditioning on the paintbox $(P_j)$, the value of $L$ is fixed, hence we are in the range of applicability of Karlin's result; see [15], Theorem 1, equation (23) and page 396. From this it is obvious that the claim holds unconditionally. □

**4. The limit distribution.** To determine the limiting behavior of $x^{\alpha_*} N_x$ as $x \to 0$ we shall connect the recursive paintbox construction to a general Crump–Mode–Jagers (CMJ) branching process [13]. The idea is to map the sizes of crumbs into a continuous time scale.

The setup for a CMJ branching process involves the random data $(\pi, \chi)$ with $\pi$ a prototypical point process on $\mathbb{R}_+$ according to which descendants are born, and $(\chi(t), t \in \mathbb{R})$ a process called *characteristic* (or a score of individual), which is nonnegative and satisfies $\chi(t) = 0$ for $t < 0$. The branching process starts at time 0 with a single progenitor which produces offsprings at epochs of $\pi$, and each descendant follows the same kind of behavior independently of the history and of the coexisting individuals. Labeling individuals in the genealogical order by integer sequences $w = j_1, \ldots, j_n$, let $\tau_w$ be the birth epoch of the generic individual. The CMJ process is defined as [13]

$$Z_t^\chi = \sum_w \chi_w(t - \tau_w),$$

which is the sum of characteristics of individuals born before $t$.

To represent the configuration of crumbs as a CMJ process we set $\tau_w = -\log \xi_w$ for crumb labeled $w$ and we define the characteristic

$$\chi_w(t) = \#\{j : -\log(\eta_{wj}/\xi_w) \leq t\}$$

to encode the configuration of solids produced by the crumb. It follows easily from the definitions that $Z_t^\chi = N_{e^{-t}}$. A key point is to apply [19], Theorem 5.4.

LEMMA 3. *As $t \to \infty$*

$$e^{-\alpha_* t} Z_t^\chi \to M \frac{\varphi(\alpha_*)}{-\alpha_* \psi'(\alpha_*)} \qquad a.s.$$

*for $M$ the terminal value of the intrinsic martingale.*

PROOF. Translated in our terms, Conditions 5.1 and 5.2 from [19] require existence of integrable, bounded, nonincreasing positive functions $h_1$ and $h_2$ such that

$$\mathbb{E}\left[\sup_t \frac{1}{h_1(t)} \sum_i X_i^{\alpha_*} 1(X_i < e^{-t})\right] < \infty$$

and

$$\mathbb{E}\left[\sup_t \frac{e^{-\alpha_* t} \#\{j : Y_j \geq e^{-t}\}}{h_2(t)}\right] < \infty.$$

These two inequalities follow from the Malthusian hypothesis with $h_1(t) = h_2(t) = e^{-\varepsilon t}$ for sufficiently small $\varepsilon > 0$. Applying [19], Theorem 5.4, we see that

$$(11) \qquad e^{-\alpha_* t} Z_t^\chi \to N \frac{\int_0^\infty e^{-\alpha_* t} \mathbb{E}[\chi(t)] \, dt}{\int_0^\infty u e^{-\alpha_* u} \mu(du)} \qquad \text{a.s.}$$

where $N$ is the terminal value of some martingale (different from the intrinsic martingale) and $\mu$ is the intensity measure of $\tau$, that is, the image of measure $\sigma$ via mapping $x \mapsto -\log x$. Since $\mathbb{E}[\chi(t)] = \int_0^1 1(x \geq e^{-t}) \nu(dx)$ by definition of intensity $\nu$, the numerator in the r.h.s. of (11) is

$$\int_0^\infty e^{-\alpha_* t} \mathbb{E}[\chi(t)] \, dt = \int_0^1 y^{\alpha_* - 1} \int_0^1 1(x \geq y) \nu(dx) \, dy$$

$$= \int_0^1 \frac{x^{\alpha_*}}{\alpha_*} \nu(dx) = \frac{\varphi(\alpha_*)}{\alpha_*}$$

by Fubini's theorem which is applicable due to the Malthusian hypothesis. Changing the variable in the numerator of the r.h.s. of (11) we see that it is equal to $-\psi'(\alpha_*)$. Hence

$$e^{-\alpha_* t} Z_t^\chi \to N \frac{\varphi(\alpha_*)}{-\alpha_* \psi'(\alpha_*)} \qquad \text{a.s.}$$

We can also apply the same result to a CMJ branching process with different characteristic $\chi'(t) = 1(t \geq 0)$, which counts individuals born before $t$. Condition 5.1 remains the same and Condition 5.2 becomes $\sup_t e^{-\alpha_* t}/h_2(t) < \infty$, so we can take $h_2(t) = e^{-\alpha_* t/2}$. Thus the Malthusian hypothesis implies that

(12) $$e^{-\alpha_* t} Z_t^{\chi'} \to N \frac{1}{-\alpha_* \psi'(\alpha_*)} \qquad \text{a.s.}$$

with the same $N$ as above.

Biggins [4] derived similar asymptotics for $Z_t^{\chi'}$ in terms of branching random walks. From [4], Theorem B and the Malthusian hypothesis

$$\frac{1}{T} \int_0^T e^{-\alpha_* t} \, dZ_t^{\chi'} \to M \frac{1}{-\psi'(\alpha*)} \qquad \text{a.s.}$$

as $T \to \infty$, where $M$ is the terminal value of the intrinsic martingale. Integration by parts and comparison with (12) show that $N = M$ a.s. $\square$

Translating the lemma back in terms of the sizes of solids we have:

COROLLARY 4. *As* $x \downarrow 0$

$$x^{\alpha_*} N_x \to \frac{\varphi(\alpha_*)}{-\alpha_* \psi'(\alpha_*)} M \qquad a.s.$$

Next, combining this corollary with Lemma 2 gives our principal asymptotic result which complements Theorem 1.

THEOREM 5. *If assumption* (1) *and the Malthusian hypothesis both hold, then*

$$K_n \sim n^{\alpha_*} \left[ \frac{\varphi(\alpha_*) \Gamma(-\alpha_*)}{\psi'(\alpha_*)} \right] M \qquad as \ n \to \infty,$$

$$K_{nr} \sim n^{\alpha_*} \left[ \frac{\varphi(\alpha_*) \Gamma(r - \alpha_*)}{-\psi'(\alpha_*) r!} \right] M \qquad as \ n \to \infty,$$

*almost surely.*

**5. Moments of $M$.** Formulas for moments of the terminal value of the intrinsic martingale involve expectations of some symmetric functions in the variables $(X_i)$. For each integer vector $\lambda = (\lambda_1, \ldots, \lambda_\ell)$ with components $\lambda_1 \geq \cdots \geq \lambda_\ell > 0$ let

$$m(\lambda) = \mathbb{E}\left[\sum_{(\mu_1,\ldots,\mu_\ell)} \sum_{i_1 < \cdots < i_\ell} X_{i_1}^{\mu_1 \alpha_*} \cdots X_{i_\ell}^{\mu_\ell \alpha_*}\right],$$

where the external sum expands over all *distinct* permutations $(\mu_1, \ldots, \mu_\ell)$ of the entries of $(\lambda_1, \ldots, \lambda_\ell)$, and the internal sum expands over all increasing $\ell$-tuples of labels of $(X_i)$. We assume for the rest of the paper that these moments exist for all integer vectors $\lambda$; this is always the case if the number of positive $X_i$'s does not exceed some constant, since $X_i \leq 1$ for all $i$.

Let $a_k = \mathbb{E}[M^k]$ ($k = 0, 1, \ldots$) be the moments of the terminal value $M$ of the intrinsic martingale (they all are finite, see [12], Proposition 4). In principle, the moments can be determined recursively from the following lemma.

LEMMA 6. *Under assumption* (1), *the Malthusian hypothesis and finiteness of moments $m(\lambda)$, the moments $a_k$ satisfy the recursion*

(13) $$a_k = \frac{k!}{1 - \psi(\alpha_* k)} \sum_{\substack{\lambda \vdash k \\ \lambda \neq (k)}} m(\lambda) \prod_j \frac{a_{\lambda_j}}{\lambda_j!} \qquad \text{for } k = 2, 3, \ldots$$

*where the initial values are $a_0 = a_1 = 1$ and the summation is over all nonincreasing positive integer sequences $\lambda = (\lambda_1, \ldots, \lambda_\ell)$ with $\lambda_1 + \cdots + \lambda_\ell = k$ and $\ell > 1$.*

PROOF. Take the $k$th power in (3) and expand the r.h.s. by the multinomial formula. Collecting all terms containing $a_k$ to the left side yields the recursion. □

## 6. Dirichlet splittings.

6.1. *Bessel bridges.* It is known that the Ewens–Pitman $(\alpha, \alpha)$ partition structure ($0 < \alpha < 1$) can be induced by a paintbox $P$ whose components are the lengths of excursions from 0 of a Bessel bridge $(B_t, t \in [0, 1])$ of dimension $2 - 2\alpha$ [7, 21, 23]. A possible recursive construction of $P$ is the following. For each $t \in [0, 1]$ define $\mathcal{G}_t = \sup\{s \leq t : B_s = 0\}$ and $\mathcal{D}_t = \inf\{s \geq t : B_s = 0\}$. Choose a random point $T$ from some distribution on $]0, 1[$, independently of $(B_t)$. The bridge $(B_t)$ decomposes into three components according as $0 \leq t \leq \mathcal{G}_T$ (bridge), $\mathcal{G}_T < t < \mathcal{D}_T$ (excursion) or $\mathcal{D}_T \leq t \leq 1$ (bridge). Given $\mathcal{G}_T$ and $\mathcal{D}_T$, the components are conditionally independent and the first and the third components are the scaled copies of $(B_t)$. It follows that the iterated division in three intervals (bridge-excursion-bridge) yields the same $P$ as the recursive construction directed by $(X_1, Y_2, X_3) \stackrel{d}{=} (\mathcal{G}_T, \mathcal{D}_T - \mathcal{G}_T, 1 - \mathcal{D}_T)$ with arbitrary distribution for $T$.

In particular, assuming $T \stackrel{d}{=} U$ for $U$ uniform $[0, 1]$, the law of $(X_1, Y_2, X_3)$ is Dirichlet with parameters $(\alpha, 1 - \alpha, \alpha)$, because this is the law of $(\mathcal{G}_U, \mathcal{D}_U - \mathcal{G}_U, 1 - \mathcal{G}_U)$, as in [21]. Computing

$$\psi(\beta) = \frac{2\alpha}{\beta + \alpha} \quad \text{and} \quad \varphi(\beta) = \frac{\Gamma(1+\alpha)\Gamma(\beta + 1 - \alpha)}{\Gamma(1-\alpha)\Gamma(\beta + 1 + \alpha)}$$

we see that in this case the Malthusian exponent is $\alpha_* = \alpha$. Applying Theorem 5 we obtain the asymptotics

$$K_n \sim \frac{\Gamma(\alpha)}{\Gamma(2\alpha)} M n^\alpha, \qquad n \to \infty,$$

which is the "$\alpha$-diversity" of $\mathcal{P}$ previously shown in [10, 22, 23] by different methods. The variable $M$ has moments

$$\mathbb{E}[M^q] = \frac{\Gamma(\alpha)\Gamma(q+1)}{[\Gamma(\alpha)/\Gamma(2\alpha)]^q \Gamma((q+1)\alpha)}, \qquad q > -1$$

and its distribution is a size-biased version of the Mittag–Leffler distribution; see [21].

Choosing any other distribution for $T$ (e.g., $T = 1/2$ a.s.) will result in different distribution for $(X_1, Y_2, X_3)$, although, by the special self-similarity properties of this $(\alpha, \alpha)$ case, the law of $P$ (up to arrangement of terms) will not alter.

We recall the original construction of the Poisson–Dirichlet paintbox from [20]; see also [22]. Let $0 < \alpha < 1$ and $\theta > -\alpha$. Take $(W_j)$ to be a sequence of independent random variables where $W_j$ has beta$(\theta + j\alpha, 1 - \alpha)$ distribution. Then the Poisson–Dirichlet $(\alpha, \theta)$ paintbox can be composed of the terms $P_j = W_1 \cdots W_{j-1}(1 - W_j)$.

6.2. *Other tripartite Dirichlet splittings.* To extend the above Bessel bridge model assume that the triple $(X_1, Y_2, X_3)$ has a Dirichlet distribution with parameters $(\gamma, \beta, \gamma)$, where $\beta, \gamma > 0$. In this case the intensity measure $\sigma$ has a density which is beta$(\gamma, \beta + \gamma)$ multiplied by 2, and $\nu$ is beta$(\beta, 2\gamma)$. Their Mellin transforms are

$$\psi(\alpha) = 2\frac{\Gamma(\beta + 2\gamma)\Gamma(\alpha + \gamma)}{\Gamma(\alpha + \beta + 2\gamma)\Gamma(\gamma)} \quad \text{and} \quad \varphi(\alpha) = \frac{\Gamma(\beta + 2\gamma)\Gamma(\alpha + \beta)}{\Gamma(\beta)\Gamma(\alpha + \beta + 2\gamma)}.$$

The recursion for moments $a_k = \mathbb{E}[M^k]$ in Lemma 6 specializes as

$$(14) \quad a_n = \sum_{k=0}^{n} \binom{n}{k} a_k a_{n-k} \frac{\Gamma(\beta + 2\gamma)\Gamma(k\alpha_* + \gamma)\Gamma((n-k)\alpha_* + \gamma)}{\Gamma(n\alpha_* + \beta + 2\gamma)\Gamma(\gamma)^2}, \qquad n \geq 2$$

with the initial values $a_0 = a_1 = 1$.

In the case when $r = \beta + \gamma$ is integer, $\psi$ is a rational function, and (2) is actually a polynomial equation in $\alpha$ of degree $r$. The case $r = 1$ covers the Bessel bridge instance of the previous section. For $r = 1, 2, \ldots$ simplification is possible by introducing variables

$$b_n = \frac{\Gamma(n\alpha_* + \gamma)}{\Gamma(\gamma)n!} a_n \tag{15}$$

for which the recursion (14) becomes

$$\sum_{k=0}^{n} b_k b_{n-k} = \frac{(n\alpha_* + \gamma)_{r\uparrow}}{(\gamma)_{r\uparrow}} b_n, \qquad n \geq 2. \tag{16}$$

Note that the same formulae also hold for $n = 0, 1$ since $b_0 = 1$ and in view of (2). This allows us to characterize the generating function

$$h(y) := \sum_{k=0}^{\infty} b_k y^k$$

as a solution to the differential equation

$$z^{1-\gamma} \frac{\mathrm{d}^r}{\mathrm{d}z^r}(h(z^{\alpha_*})z^{\gamma+r-1}) = (\gamma)_{r\uparrow} h(z^{\alpha_*})^2. \tag{17}$$

In the variable $y = z^{\alpha_*}$ this equation is a nonlinear differential equation with polynomial coefficients.

For instance, when $r = 2$, (2) becomes $(\alpha_* + \gamma)(\alpha_* + \gamma + 1) = 2\gamma(\gamma + 1)$ and after some manipulations we obtain

$$\alpha_*^2 y^2 h''(y)/(\gamma(\gamma+1)) + yh'(y) + h(y) = h^2(y).$$

We did not succeed in solving the equation in terms of some known special functions for $r \geq 2$. We can, nevertheless, show that this partition structure is of novel type:

LEMMA 7. *For no $r = 2, 3, \ldots$ and no $\gamma \in\, ]0, r[$ does the recursive partition structure obtained by the recursive tripartite Dirichlet splitting with parameters $(\gamma, r - \gamma, \gamma)$ belong to the Ewens–Pitman two-parameter family of partition structures.*

PROOF. The statement follows by computing probabilities $p(n)$ for $n$ balls painted the same color. Indeed, in the $(\alpha, \theta)$-model this probability is [22]

$$p_{\alpha,\theta}(n) = \frac{(1-\alpha)_{(n-1)\uparrow}}{(1+\theta)_{(n-1)\uparrow}},$$

and in our model it is

$$p(n) = \frac{(r-\gamma)_{n\uparrow}}{(r+\gamma)_{n\uparrow} - 2(\gamma)_{n\uparrow}}.$$

Assuming the coincidence for some value of parameters $(\alpha, \theta)$ we must have $\alpha = \alpha_*$ and $p(n) = p_{\alpha,\theta}(n)$ for all $n$. Analyzing the behavior of these probabilities as $n \to \infty$ we find out that

$$p_{\alpha,\theta}(n) \sim \frac{\Gamma(1+\theta)}{\Gamma(1-\alpha)} n^{-\alpha-\theta} \quad \text{and} \quad p(n) \sim \frac{\Gamma(r+\gamma)}{\Gamma(r-\gamma)} n^{-2\gamma},$$

whence $\theta = 2\gamma - \alpha$. Substituting this value in the equation $p(2) = p_{\alpha,\theta}(2)$ we see that $\alpha = \gamma - \frac{r^2-r}{2r-\gamma}$. Comparing again the rates of decrease of $p_{\alpha,\theta}(n)$ and $p(n)$ for these particular values of $\alpha, \theta$ we get

$$\frac{\Gamma(1+(r^2-r)/(2r-\gamma)+\gamma)}{\Gamma(1+(r^2-r)/(2r-\gamma)-\gamma)} = \frac{\Gamma(r+\gamma)}{\Gamma(r-\gamma)}.$$

Since the r.h.s. increases as a function of $r > \gamma$ and the l.h.s. is the same function evaluated at a different point, necessarily $1 + \frac{r^2-r}{2r-\gamma} = r$. But this can happen only for $r = 1$ or $r = \gamma$. In the latter case the parameter $\gamma = r$ is not admissible, so the coincidence happens only for $r = 1$. $\square$

6.3. *Multiple splittings and $(\alpha, \alpha/d)$ partitions.* Now suppose the splitting procedure produces $d+1$ crumbs $(d \geq 1)$ and one solid at each step. Suppose the joint distribution of the crumb sizes and the solid size relative to their parent crumb is the Dirichlet distribution with parameters $(\gamma, \ldots, \gamma, \beta)$ where $\gamma$'s correspond to $d+1$ crumbs and $\beta$ corresponds to the solid. Mellin transforms of the intensity measures are

$$\psi(\alpha) = (d+1) \frac{\Gamma(\beta + (d+1)\gamma)\Gamma(\alpha+\gamma)}{\Gamma(\alpha+\beta+(d+1)\gamma)\Gamma(\gamma)}$$

and

$$\varphi(\alpha) = \frac{\Gamma(\beta+(d+1)\gamma)\Gamma(\alpha+\beta)}{\Gamma(\beta)\Gamma(\alpha+\beta+(d+1)\gamma)}.$$

The $d = 1$ case was considered in the preceding section. Similarly to the above, explicit computations are only possible when $\beta + d\gamma = r$ is integer.

The simplest case $r = 1$ leads to some exactly solvable recursion for moments of the terminal value $M$ of the intrinsic martingale. We consider this case in more detail. For $r = 1$, (2) becomes $(d+1)\gamma/(\alpha+\gamma) = 1$ with the solution $\alpha_* = d\gamma$. The recursion for moments $a_n$ of $M$ is easier to write down in new variables $b_n$ defined by (15):

(18) $$\sum b_{\lambda_1} \ldots b_{\lambda_{d+1}} = (nd+1)b_n, \qquad n \geq 2,$$

where the sum is taken over all nonnegative integer vectors $(\lambda_1, \ldots, \lambda_{d+1})$ with $\sum_j \lambda_j = n$. This leads to the differential equation $dyh'(y) + h(y) = h(y)^{d+1}$ for the generating function $h(y) = \sum_{n=0}^{\infty} b_n y^n$. Solving this equation we obtain

$$\mathbb{E}[M^q] = \frac{d^q \Gamma(\alpha_* + \alpha_*/d)^q \Gamma(1/d + q)}{\Gamma(\alpha_*/d)^{q-1} \Gamma(q\alpha_* + \alpha_*/d) \Gamma(1/d)}.$$

We recognize these as the moments of the limit distribution of $\frac{d\Gamma(\alpha_* + \alpha_*/d)}{\Gamma(\alpha_*/d)} n^{-\alpha_*} K_n$, where $K_n$ is the number of blocks in the Ewens–Pitman partition structure with parameters $(\alpha_*, \alpha_*/d)$ [20, 22] restricted to the first $n$ balls. The limit has density proportional to $x^{1/d} f_{\alpha_*}(x)$, where $f_\alpha$ is the density of the Mittag–Leffler distribution with parameter $\alpha$. The following proposition shows that the partition structures coincide.

PROPOSITION 8. *The exchangeable partition obtained by a splitting scheme with $d+1$ crumbs and one solid whose joint distribution is* Dirichlet $(\underbrace{\alpha/d, \ldots, \alpha/d}_{d+1}, 1-\alpha)$ $(0 < \alpha < 1)$ *coincides with the Ewens–Pitman $(\alpha, \alpha/d)$ partition.*

In the proof we use a mapping $q$ which sends a collection $B$ of $kd+1$ elements with unit sum to a random collection of $(k+1)d+1$ elements with unit sum. This mapping is defined for $\alpha \in\, ]0, 1[$ as follows:

(1) choose an element from the collection $B$ by a size-biased pick;
(2) replace the chosen element $Z$ by $d+2$ elements $(YZ, X_1Z, \ldots, X_{d+1}Z)$ where $(Y, X_1, \ldots, X_{d+1})$ is an independent of $B$ $Dirichlet(1-\alpha, \alpha/d, \ldots, \alpha/d)$ random vector;
(3) remove the element $YZ$ from the collection, divide all elements by $1-YZ$ so that they sum to 1, and let $q(B)$ be the rescaled collection.

LEMMA 9. *Let $B$ be a collection of $kd + 1$ random variables whose joint distribution is* Dirichlet$(\alpha/d, \ldots, \alpha/d)$. *Then $q(B)$ is a collection of $(k+1)d+1$ random variables with joint distribution* Dirichlet$(\alpha/d, \ldots, \alpha/d)$. *Moreover, the size of the discarded element is a* beta$(1 - \alpha, (k + 1 + 1/d)\alpha)$ *random variable independent of $q(B)$.*

PROOF. After the first step, the conditional distribution of elements in $B$ given that the size-biased pick has index $i$ is the Dirichlet distribution with one parameter $\alpha/d + 1$ for element $i$ and other parameters $\alpha/d$. We relabel the elements so that the chosen element is the first one. After the second step, the elements in the collection have the Dirichlet distribution with the first parameter $1 - \alpha$ and other $(k+1)d + 1$ parameters $\alpha/d$. This can be easily verified by a moment calculation, exploiting the independence of $(Y, X_1, \ldots, X_{d+1})$ and $B$ and the fact that the Dirichlet parameters of $(Y, X_1, \ldots, X_{d+1})$ sum to the Dirichlet parameter of the replaced element. The statement of the lemma now follows from [14], Chapter 40. $\square$

PROOF OF PROPOSITION 8. The Poisson–Dirichlet paintbox can be arranged in a sequence by the stick-breaking procedure described in the end of Section 6.1. We show that the terms of the paintbox in our model can be also arranged in such sequence, as follows. Since each crumb produces exactly one solid in the model in focus, there is a one-to-one correspondence between solids and their parent crumbs. Let the first solid $\eta(1)$ in the arrangement be a child of the progenitor crumb $\xi_\varnothing$ and let $A_1 = \{\xi_1, \ldots, \xi_{d+1}\}$ be the offspring crumbs of $\xi_\varnothing$. Inductively, at time $k$ let the first $k$ solids have been arranged as $\eta(1), \ldots, \eta(k)$ and let $A_k$ be some collection of crumbs. The next solid to be added to the sequence is chosen in the following way. Select a crumb $\xi_w$ by a size-biased pick from all crumbs in the collection $A_k$, let the next element $\eta(k+1)$ added to $P$ be the solid child of this $\xi_w$, and further replace $\xi_w$ in the collection $A_k$ by the offspring crumbs of $\xi_w$, thus constructing $A_{k+1} := (A_k \setminus \{\xi_w\}) \cup \{\xi_{w,1}, \ldots, \xi_{w,(d+1)}\}$. Proceed by induction to arrange all solids in sequence.

Now let us check that the sequence of solids $P = (\eta(i))$ has the same distribution as the lengths produced by the stick-breaking procedure described at the end of Section 6.1. At the first step, the law of $\eta(1)$ is the marginal distribution of the Dirichlet distribution which is beta$(1 - \alpha, (1 + 1/d)\alpha)$. For $k = 1, 2, \ldots$ introduce the scaled collections of crumbs

$$B_k = \left\{ \frac{\xi}{|A_k|} : \xi \in A_k \right\}, \qquad \text{where } |A_k| = \sum_{\xi \in A_k} \xi = 1 - \eta(1) - \cdots - \eta(k).$$

Then $B_1$ has Dirichlet$(\alpha/d, \ldots, \alpha/d)$ distribution by [14], Chapter 40, and $B_{k+1} = q(B_k)$ for all $k$. Using Lemma 9, we check by induction that $B_k$ is a collection of $dk + 1$ elements whose joint distribution is Dirichlet with all parameters $\alpha/d$, and $\eta(k)$ has beta$(1 - \alpha, (k + 1/d)\alpha)$ distribution and is independent of $\eta(1), \ldots, \eta(k-1)$ for all $k$. Taking $W_k = 1 - \eta(k)$ yields the desired decomposition. $\square$

**7. Further subdivision of solids.** Suppose we have some recursive paintbox construction with the Mellin transforms $\psi_0$ and $\varphi_0$ of the intensity measures for **X** and **Y**. A refined paintbox construction can be produced by a further independent subdivision of each solid according to some sequence $\tilde{P} = (\tilde{Y}_k)$ of nonnegative random variables with $\sum \tilde{Y}_k = 1$. This is equivalent to replacing $(Y_j)$ in the original construction by an array $(Y_j \tilde{Y}_k)$ (arranged in a sequence). By independence, the Mellin transform of a new intensity measure for $(Y_j \tilde{Y}_k)$ in the refined process is the product

$$\varphi(\alpha) = \varphi_0(\alpha)\tilde{\varphi}(\alpha), \qquad \tilde{\varphi}(\alpha) = \mathbb{E}\left[\sum_k \tilde{Y}_k^\alpha\right].$$

If the expected number of nonzero $\tilde{Y}_k$'s is finite, then this new construction satisfies the Malthusian hypothesis once the original construction satisfied it. If an infinite number of positive $\tilde{Y}_k$'s is possible, we should also require $\tilde{\varphi}(\alpha_* - \varepsilon) < \infty$ to keep with the Malthusian hypothesis.

One example where a similar additional subdivision of solids was used is a representation of the Poisson–Dirichlet $(\alpha, \theta)$ paintbox [although it does not fit exactly in our scheme since the expected number of nonzero $X_i$'s is 1, violating (1)]. The recipe is the following [7, 9, 24]: divide the unit interval by points of a stick-breaking process with $W_i$ i.i.d. beta$(\theta, 1)$ and then organize on each subinterval an independent subdivision by zeroes of a Bessel process of dimension $2 - 2\alpha$. Here $\alpha_* = 0$ [due to a violation of (1)] and $\tilde{\varphi}(\alpha_*) = \infty$.

When the Malthusian hypothesis still holds for the refined process, it has some common features with the original one. For instance, the Malthusian exponent remains the same, and the limit of $n^{-\alpha_*} K_n$ changes only by a constant factor $\tilde{\varphi}(\alpha_*)$. However, other characteristics of the paintbox, such as probabilities $p(n)$ that $n$ balls are painted in the same color, change significantly once any subdivision is made.

**Acknowledgment.** Helpful comments of a referee are gratefully acknowledged.

DEPARTMENT OF MATHEMATICS
UTRECHT UNIVERSITY
POSTBUS 80010
3508 TA UTRECHT
THE NETHERLANDS
E-MAIL: gnedin@math.uu.nl
       yakubovich@math.uu.nl